\documentclass[10pt]{article}
\usepackage{amsmath,amssymb,amsthm,mathrsfs,dsfont}
\usepackage{needspace}

\usepackage[margin=3cm]{geometry} 
\usepackage{titlesec,hyperref}
\usepackage{amssymb}
\newcommand{\Rmnum}[1]{\uppercase\expandafter{\romannumeral #1}}  
\usepackage{color}

\usepackage{fancyhdr}
\titleformat{\subsection}{\it}{\thesubsection.\enspace}{1pt}{}

\newtheorem{theo}{Theorem}[section]
\newtheorem{lemm}[theo]{Lemma}
\newtheorem{defi}[theo]{Definition}

\newtheorem{prop}[theo]{Proposition}
\newtheorem{rema}[theo]{Remark}
\numberwithin{equation}{section}

\allowdisplaybreaks 

\newcommand\lm{{\lesssim}}

\begin{document}
\title{Global existence and sharp decay rates of weak solutions for the 2D inviscid Oldroyd-B model with
fractional stress tensor diffusion
\hspace{-4mm}
}

\author{Chen $\mbox{Liang}^1$ \footnote{Email: liangch89@mail2.sysu.edu.cn}\quad and\quad
	Zhaonan $\mbox{Luo}^1$\footnote{E-mail: luozhn7@mail.sysu.edu.cn} \\
	$^1\mbox{School}$ of Science,\\ Shenzhen Campus of Sun Yat-sen University, Shenzhen 518107, China}
\date{}
\maketitle
\hrule

\begin{abstract}
This paper is devoted to the global existence and optimal decay rates of weak solutions to a two-dimensional inviscid Oldroyd-B model with fractional stress diffusion. For $\frac{1}{2}<\beta<1$, by virtue of energy estimates and compactness arguments, we establish the existence of global weak solutions with $(u,\tau)\in H^{1}\times (H^{1}\cap L^{\infty})$. Moreover, based on the Littlewood-Paley decomposition theory and the improved Fourier splitting method, we derive the optimal decay rates. In particular, by constructing a suitable low-frequency energy functional and exploiting the dissipation mechanism of the coupled system, we obtain the optimal $L^\infty$ decay rate for the stress tensor.\\
\vspace*{5pt}

\vspace*{5pt}
\end{abstract}

\vspace*{-20pt}

	\section{Introduction}

We consider the following two-dimensional incompressible Oldroyd-B system with only fractional stress diffusion(\cite{2015Elgindi,1958Non}):
	\begin{align}\label{eq1}
\left\{\begin{array}{l}
			\partial_tu + u\cdot\nabla u+\nabla {\rm P} = {\rm div}~\tau,~~~~{\rm div}~u=0,\\[1ex]
			\partial_t\tau + u\cdot\nabla\tau+a\tau+Q(\nabla u,\tau)+(-\Delta)^{\beta}\tau=D(u),\\[1ex]
			u|_{t=0}=u_0,~~\tau|_{t=0}=\tau_0. \\[1ex]
		\end{array}\right.
	\end{align}
	Here, $u(t,x)$ denotes the velocity field,
$\tau(t,x)$ is the symmetric extra-stress tensor, and $\rm P$ denotes the pressure. The bilinear form $Q$ is given by
	$$Q(\nabla u, \tau)=\tau \Omega(u)-\Omega(u)\tau-b(D(u)\tau+\tau D(u)),~~~~~~~b\in[-1, 1],$$
	with vorticity tensor
	$$
	\Omega(u)=\frac {\nabla u-(\nabla u)^T} {2},
	$$
	and deformation tensor
	$$
	D(u)=\frac {\nabla u+(\nabla u)^T} {2}.
	$$
    $0\le a\le 1$ and $\frac{1}{2}\le\beta\le1$ are parameters. Fractional Laplacian operator $(-\Delta)^{\beta}$ is defined through the Fourier transform by
	$$
	(-\Delta)^{\beta}\tau=\mathscr{F}^{-1}(|\xi|^{2\beta}\widehat{\tau}),\quad\quad \widehat{\tau}(\xi,t)=\mathscr{F}\tau(\xi,t)=\int_{\mathbb{R}^2}e^{-ix\cdot\xi}\tau(x,t) dx.
	$$
When $\beta=1$, system (\ref{eq1}) is the classical  Oldroyd-B model. The Oldroyd-B model is widely used to describe the motion of viscoelastic fluids. For more details on the classical Oldroyd-B model, we refer to \cite{1958Non,Oldroyd1950On}.

For the case $\beta=1$, the global existence and regularity of solutions to the inviscid Oldroyd-B system (\ref{eq1}) have been extensively studied in recent years. T. M. Elgindi and F. Rousset \cite{2015Elgindi} firstly proved the global existence for (\ref{eq1}) in the Sobolev space $H^s$ with $s>2$ provided that $\|(u_0,\tau_0)\|_{H^1}+\|(\mathrm{curl}\ u_0,\tau_0)\|_{B_{\infty,1}^{0}}$ is sufficiently small. Furthermore, T. M. Elgindi and J. Liu \cite{2015Elgindi1} established the global existence of small solutions in three dimensions.  When $a=0$, W. Deng, Z. Luo and Z. Yin \cite{deng2023} established the global existence of weak solutions for system (\ref{eq1}) with small initial data in $H^1$. For $a>0$, X. Cheng, Z. Luo, Z. Yang and C. Yuan \cite{CLYY} proved the global well-posedness in two dimensions under the smallness assumption of $\|(u_0,\tau_0)\|_{L^2}+\|(\nabla u_0,\tau_0)\|_{B_{\infty,1}^{0}}$. Moreover, they established the global existence of solutions in both two and three dimensions for $a\in[0,1]$ under suitable smallness assumptions. Furthermore, they obtained the sharp uniform-in-time vanishing damping rate $a^{\frac{d}{4}}$. More recently, the same authors \cite{ChengLuoYangYuan2026} established convergence and optimal convergence rates for a
low-regularity numerical scheme for the two-dimensional inviscid Oldroyd-B system with stress diffusion.

 When $a=0$ and $\frac12\leq\beta\leq1$, P. Constantin, J. Wu, J. Zhao and Y. Zhu \cite{P.Constantin} proved global well-posedness to (\ref{eq1}) with small initial data in Sobolev space $H^s$ for $s>1+\frac{d}{2}$. Moreover, J. Wu and J. Zhao \cite{Wu} obtained global well-posedness of small initial data in the critical $L^2$ framework.
 Z. Chen, L. Liu, D. Qin and W. Ye \cite{Chenzhi} further extended these results to the general $L^p$ framework for $\beta=1$.
P. Wang, J. Wu, X. Xu and Y. Zhong \cite{Wu1} studied the sharp decay rates for system (\ref{eq1}). For sufficiently small initial data in $L^1\cap H^{4+\frac{6}{\beta}}$, they obtained the following decay estimates for $\frac12\leq\beta<1$ in two dimensions:
\begin{align*}
&\|u(t)\|_{L^{2}}\lesssim\varepsilon(1+t)^{-\frac{1}{2\beta}},\quad\|u(t)\|_{L^{\infty}}\lesssim \varepsilon(1+t)^{-\frac{1}{\beta}},\quad\|\nabla u(t )\|_{L^{2}}\lesssim\varepsilon(1+t)^{-\frac{1}{\beta}},\\
&\|\nabla u(t)\|_{L^{\infty}}\lesssim\varepsilon(1+t)^{-\frac{3}{2\beta}},\quad\|\mathbb{P}\nabla\cdot\tau(t)\|_{L^ {2}}\lesssim\varepsilon(1+t)^{-\frac{1}{\beta}}.
\end{align*}
This result was further extended by Z. Luo, W. Luo and Z. Yin \cite{LLY}, who established the optimal decay rates in $H^s$ with $s>2$ for the two-dimensional system with
$\frac{1}{2}\leq\beta\leq1$.
Notably, their result includes the classical case $\beta=1$.  Recently, C. Liang, Z. Luo and Z. Yin \cite{liang2025} established the optimal decay rates for system (\ref{eq1}) with $\frac12<\beta<1$ in the critical Sobolev space $H^{2\beta}$. Moreover, they obtained the sharp uniform-in-time vanishing damping rate $a^{\frac{1}{2\beta}}$ for $\operatorname{tr}\tau^a$.

Motivated by \cite{deng2023}, we aim to study the global existence and large-time behavior of weak solutions to the system (\ref{eq1}) with 
$\frac{1}{2}<\beta<1$. We first introduce a definition of weak solutions to the system (\ref{eq1}).
\begin{defi}
Let $\phi\in C_c^\infty([0,T)\times\mathbb R^2;\mathbb R^2)$
with $\mathrm{div}\phi=0$ and
$\psi\in C_c^\infty([0,T)\times\mathbb R^2;\mathbb R^{2\times2})$.
We say that $(u^a,\tau^a)$ is a weak solution to system
(\ref{eq1}) if the following conditions hold:

(1): For any $T>0$, $(u^a,\tau^a)$ belongs to $L_T^\infty(H^1)\times L_T^\infty(H^1\cap L^\infty)$, and $\operatorname{div}u^a=0$ in the sense of distributions.

(2): For any $\phi$ and $\psi$, there holds
\begin{align*}
&\int_0^T\!\int_{\mathbb R^2}
\big(u^a\cdot\partial_s\phi+(u^a\otimes u^a):\nabla\phi-\tau^a:\nabla\phi\big)\,dx\,ds
+\int_{\mathbb R^2}u_0\cdot\phi(0)\,dx=0,
\end{align*}
and
\begin{align*}
&\int_0^T\!\int_{\mathbb R^2}\big(\tau^a:\partial_s\psi+(u^a\otimes\tau^a):\nabla\psi-a\tau^a:\psi-\tau^a:\Lambda^{2\beta}\psi
\nonumber\\
&\hspace{22mm}
-Q(\nabla u^a,\tau^a):\psi+D(u^a):\psi\big)\,dx\,ds
+\int_{\mathbb R^2}\tau_0:\psi(0)\,dx=0.
\end{align*}
\end{defi}
We first establish global regularity in Sobolev spaces for any $a\in[0,1]$. In 2D, the divergence-free condition $\mathrm{div}u^a=0$ yields the following cancellation:
\begin{align*}
\langle u^a\cdot\nabla u^a,\Delta u^a\rangle=0.
\end{align*}
This cancellation plays a crucial role in closing the $H^1$ energy estimates. On the other hand, the main difficulty in closing the $H^1$ energy estimates arises from the nonlinear term $Q(\nabla u^a,\tau^a)$, which can be estimated as follows:
\begin{align*}
\langle\Lambda Q(\nabla u^a,\tau^a),\Lambda\tau^a\rangle
\lesssim\|\Lambda^{1+\beta}\tau^a\|_{L^2}\|\Lambda^{2-\beta}u^a\|_{L^2}
\left(\|\tau^a\|_{L^\infty}+\|\Lambda\tau^a\|_{L^2}\right).
\end{align*}
Therefore, an additional $L^\infty$ estimate for $\tau^a$ is required to close the $H^1$ energy estimates. Combining the energy estimates with a continuity argument, we obtain global a priori estimates. Then, by a compactness argument, we establish the existence of global weak solutions to system \eqref{eq1}.  Our first main result is stated as follows.
\begin{theo}\label{1theo1}
Let $a\in[0,1]$ and $\frac12<\beta<1$.
Assume that $u_0\in H^1$ and $\tau_0\in H^1\cap L^\infty$ where $\operatorname{div}u_0=0$ and $\tau_0$ is symmetric. Then there exists a small constant $\delta $ such that if
\begin{align*}
\|(u_{0},\tau_{0})\|_{H^{1}}+\|\tau_0\|_{L^{\infty}}\leq\delta,
\end{align*}
then the system (\ref{eq1}) admits a global weak solution $(u^a,\tau^a)\in L^{\infty}([0,\infty),H^{1}\times (H^1\cap L^{\infty}))$. Moreover, we obtain that for all $t>0$, there holds
\begin{align*}
\quad\frac{d}{dt}\left(\|(u^a,\tau^a)\|_{H^{1}}^{2}+2k\langle-\Lambda^{2-2\beta}\nabla u^a,\tau^a\rangle\right)+\frac{k}{8}\|\Lambda^{2-\beta} u^a\|_{L^{2}}^{2}+\frac{1}{2}\|\Lambda^{\beta}\tau^a\|_{H^{1}}^{2}\leq 0,
\end{align*}
where $k$ is a sufficiently small constant.
\end{theo}

Next, we investigate the optimal decay rates of weak solutions to system \eqref{eq1} when $a=0$. By virtue of the improved Fourier splitting method and time-weighted energy estimates, we obtain the optimal decay rates in $H^1$. The main difficulty here is to establish the optimal decay rate for $\|\tau^0\|_{L^\infty}$. By Littlewood-Paley decomposition, we write
\[
\tau^0=\Delta_{-1}\tau^0+\sum_{j\ge0}\dot\Delta_j\tau^0.
\]
We observe that $\dot{\Delta}_j\tau^0$ exhibits a damping effect for $j\geq0$. Therefore, the main difficulty is to establish the decay estimate for $\|\Delta_{-1}\tau^0\|_{L^{\infty}}$. To overcome this difficulty, we introduce the following energy functional:
\[
 E_{-1}(t)=|\widehat{\Delta_{-1}u^0}|^2
+|\widehat{\Delta_{-1}\tau^0}|^2
+2k|\xi|^{2\beta-2}\operatorname{Re}
\left(\widehat{\Delta_{-1}u^0}\cdot
\overline{\widehat{\Delta_{-1}\operatorname{div}\tau^0}}\right),
\]
where $k>0$ sufficiently small. Since $\widehat{\Delta_{-1}f}$ is compactly supported in the Fourier space, for sufficiently small $k>0$, $E_{-1}(t)$ is equivalent to $|\widehat{\Delta_{-1}u^0}|^2+|\widehat{\Delta_{-1}\tau^0}|^2$. Then we derive 
\begin{align*}
\frac{d}{dt}E_{-1}^{\frac{1}{2}}+c|\xi|^{2\beta}E_{-1}^{\frac{1}{2}}
\lesssim 
|\widehat{\Delta_{-1}F}|+|\widehat{\Delta_{-1}G}|.
\end{align*}
Here, $F$ and $G$ are defined in the proof of Proposition \ref{3prop2}. This inequality reveals a parabolic dissipation mechanism for the low-frequency part of $(u^0,\tau^0)$. By Fourier's inversion formula and the above dissipation mechanism, we obtain
$$
\begin{aligned}
\|\Delta_{-1}\tau^0(t)\|_{L^\infty}
&\lesssim \int_{|\xi|\leq\frac43}
|\widehat{\Delta_{-1}\tau^0}(t,\xi)|\,d\xi\lesssim (1+t)^{-\frac1\beta}.
\end{aligned}
$$
The details are given in the proof of Proposition \ref{3prop5}. We believe our method may also be applied to other hyperbolic-parabolic model in low regularity function spaces. Based on the above analysis, we state our second main result as follows.
\begin{theo}\label{1theo2}
Let $a=0$. Let $(u^0,\tau^0)$ be a weak solution of (\ref{eq1}) with the initial data  $\left(u_{0}, \tau_{0}\right)$ under the condition in Theorem \ref{1theo1}. If additionally $\left(u_{0}, \tau_{0}\right) \in \dot{B}_{2, \infty}^{-1}$, then there exists  $C>0$  such that for every $t>0$ and $s_1\in[0,1]$,
\begin{align*}
\left\|\Lambda^{s_1}(u^0, \tau^0)\right\|_{L^{2}} \leq C(1+t)^{-\frac{s_1+1}{2 \beta}}.
\end{align*}
Moreover, we have
\begin{align*}
\left\|\tau^0\right\|_{L^{\infty}} \leq C(1+t)^{-\frac{1}{\beta}}.
\end{align*}
\end{theo}
\begin{rema}
From the proof of Theorem \ref{1theo2}, we also get
\begin{align*}
\|\Delta_{-1}u^{0}\|_{L^{\infty}}\le C(1+t)^{-\frac{1}{\beta}}.
\end{align*}
\end{rema}
\begin{rema}
All decay rates in Theorem~\ref{1theo2} are optimal (see also \cite[Theorem~1.4]{liang2025}).
In particular, for $d=2$ and $\frac12<\beta<1$, Theorem~\ref{1theo2} yields the same decay rates for
$\|u(t)\|_{L^2}$, $\|\nabla u(t)\|_{L^2}$, and $\|\mathbb P\operatorname{div}\tau(t)\|_{L^2}$ for weak solutions as those obtained in \cite[Theorem~1.3]{Wu1} for strong solutions.
\end{rema}
Several questions remain open. For the endpoint case $\beta=\frac12$, the global existence of weak solutions to system \eqref{eq1} remains open, since the condition $\frac1\beta<p<2$ required in \eqref{220} cannot be satisfied. Moreover, the uniqueness of the weak solutions obtained in Theorem \ref{1theo1} is still unknown.

{\rm\textbf{Structure of the paper:}}~~The remainder of this paper is organized as follows. In Section 2, we give some preliminaries which will be used in the sequel. In Section 3, we establish the global existence of weak solutions to system (\ref{eq1}) for $a\in[0,1]$. In Section 4, we derive the optimal decay rates for the case $a=0$.

	\section{Preliminaries}
\par
	In this section, we introduce some notations and useful lemmas which will be used in the sequel.
	
	We agree that $f\lm g$ represents $f\leq Cg$ with a constant $C>0$. The symbol $\widehat{f}=\mathscr{F}(f)$ stands for the Fourier transform of $f$. Denote $\mathscr{F}^{-1}(f)$ the inverse Fourier transform of $f$.
	
	The Littlewood-Paley decomposition theory and Besov spaces are given as follows:
	\begin{lemm}\cite{Bahouri2011}\label{lemma1}
		Let $\mathscr{C}$ be the annulus $\{\xi\in\mathbb{R}^2:\frac 3 4\leq|\xi|\leq\frac 8 3\}$. There exist  radial functions $\chi$ and $\varphi$, valued in the interval $[0,1]$, belonging respectively to $D(B(0,4/3))$ and $\mathscr{D}(\mathscr{C})$, and such that
		$$ \forall\xi\in\mathbb{R}^2,\ \chi(\xi)+\sum_{j\geq 0}\varphi(2^{-j}\xi)=1, $$
		$$ \forall\xi\in\mathbb{R}^2\backslash\{0\},\ \sum_{j\in\mathbb{Z}}\varphi(2^{-j}\xi)=1,~~~ $$
		$$ |j-j'|\geq 2\Rightarrow\mathrm{Supp}\ \varphi(2^{-j}\cdot)\cap \mathrm{Supp}\ \varphi(2^{-j'}\cdot)=\emptyset, $$
		Moreover, we have
		$$ \forall\xi\in\mathbb{R}^2\backslash\{0\},\ \frac 1 2\leq\sum_{j\in\mathbb{Z}}\varphi^2(2^{-j}\xi)\leq 1.~~ $$
	\end{lemm}
			
		Let $u$ be a tempered distribution in $\mathcal{S}_{h}'(\mathbb{R}^2)$. For all $j\in\mathbb{Z}$, define
		The homogeneous operators are defined by
		$$\dot{\Delta}_j u=\mathscr{F}^{-1}(\varphi(2^{-j}\cdot)\mathscr{F}u).$$
Then the Littlewood-Paley decomposition is given as follows:
$$
u=\sum_{j\in\mathbb{Z}}\dot{\Delta}_j u\quad  \mathrm{in}\quad  \mathcal{S}_{h}'(\mathbb{R}^2).
$$		
		Let $s\in\mathbb{R}$ and $(p,r)\in[1,\infty]^2$. The homogeneous Besov space $\dot{B}^s_{p,r}$ is given as follows
	$$ \dot{B}^s_{p,r}=\left\{u\in \mathcal{S}_{h}'(\mathbb{R}^2):\|u\|_{\dot{B}^s_{p,r}}=\Big\|(2^{js}\|\dot{\Delta}_j u\|_{L^p})_j \Big\|_{l^r(\mathbb{Z})}<\infty\right\}.$$

Let
		$${\Delta}_{-1} u=\mathscr{F}^{-1}(\chi(\cdot)\mathscr{F}u)$$
        Then the Littlewood-Paley decomposition is given as follows:
$$
u=\Delta_{-1}u+\sum_{j\ge 0}{\dot\Delta}_j u\quad  \mathrm{in}\quad  \mathcal{S}'(\mathbb{R}^2).
$$		
\begin{lemm}
\label{lem:heat-smoothing}\cite{Bahouri2011}
Let $\mathcal{C}$ be an annulus. Positive constants $c$ and $C$ exist such that for any $p$ in $[1,\infty]$ and any couple $(t,\lambda)$ of positive real numbers, we have
\[
\operatorname{Supp}\widehat{u}\subset \lambda \mathcal{C}
\;\Longrightarrow\;
\| e^{-t(-\Delta)^{\beta}}u\|_{L^{p}}\leq C e^{-c t\lambda^{2\beta}}\| u\|_{L^{p}}.
\]
\end{lemm}
The following commutator estimate and product estimate are useful to control convective terms and $Q(\nabla u,\tau)$.
\begin{lemm}\label{lemma4}\cite{kato1}
Assume that $s>0$, $p, p_1, p_4 \in (1,\infty)$ and $\frac{1}{p}=\frac{1}{p_1}+\frac{1}{p_2}=\frac{1}{p_3}+\frac{1}{p_4}$, then we obtain
\begin{align*}
\|[\Lambda^{s},f]g\|_{L^{p}}\leq C\left(\|\Lambda^{s}f\|_{L^{p_{1}}}\|g\|_{L^{p_{2} }}+\|\nabla f\|_{L^{p_{3}}}\|\Lambda^{s-1}g\|_{L^{p_{4}}}\right).
\end{align*}
\end{lemm}
\begin{lemm}\label{lemma5}\cite{kato1}
Assume that $s>0$, $p, p_2, p_4 \in (1,\infty)$ and $\frac{1}{p}=\frac{1}{p_1}+\frac{1}{p_2}=\frac{1}{p_3}+\frac{1}{p_4}$, then we obtain
\begin{align*}
\|\Lambda^{s}\left(fg\right)\|_{L^{p}} \leq C\left(\|f \|_{L^{p_{1}}}\|\Lambda^{s}g\|_{L^{p_{2}}} + \|g\|_{L^{p_ {3}}}\|\Lambda^{s}f\|_{L^{p_{4}}}\right).
\end{align*}
\end{lemm}
The following fractional heat semigroup estimates will be used to bound the stress tensor in $L^\infty$.

\section{Global existence of weak solutions}
In this section, we establish the global existence of weak solutions to system \eqref{eq1} for $a\in[0,1]$ via uniform a priori estimates and compactness arguments. We first establish the a priori estimates for system \eqref{eq1}.

 \textbf{Proof of the $H^1$ energy estimate:}\\
 Let $(u^a,\tau^a)$ be a local smooth solution of equations (\ref{eq1}). We assume a priori that
 \begin{align}\label{210}
\|(u^a,\tau^a)\|_{H^1}+\|\tau^a\|_{L^\infty}\le\delta.  
 \end{align}
Taking the $L^2$ inner product to system (\ref{eq1}) with $(u^a,\tau^a)$ and integrating by parts, we get
\begin{align*}
\quad\frac{1}{2}\frac{d}{dt}\|(u^a,\tau^a)\|_{L^{2}}^{2}+\|\Lambda^{\beta}\tau^a\|_{L^{2}}^{2}+a\|\tau^a\|_{L^{2}}^{2}=-\langle Q(\nabla u^a,\tau^a), \tau^a\rangle,
\end{align*}
where we have used the cancellation
\begin{align*}
\langle u^a,\mathrm{div}\tau^a\rangle+\langle Du^a,\tau^a\rangle=0.
\end{align*}
Since
\begin{align*}
\langle Q(\nabla u^a,\tau^a), \tau^a\rangle\lesssim \ \|\nabla u^a\|_{L^{\frac{2}{\beta}}}\|\tau^a\|_{L^{\frac{2}{1-\beta}}}\|\tau^a\|_{L^{\infty}}
\lesssim \ \|\Lambda^{2-\beta} u^a\|_{L^2}\|\Lambda^{\beta}\tau^a\|_{L^2}\|\tau^a\|_{L^{\infty}},
\end{align*}
then we deduce that
\begin{align}\label{211}
\quad\frac{1}{2}\frac{d}{dt}\|(u^a,\tau^a)\|_{L^{2}}^{2}+\|\Lambda^{\beta}\tau^a\|_{L^{2}}^{2}+a\|\tau^a\|_{L^{2}}^{2}
\lesssim\|\Lambda^{2-\beta} u^a\|_{L^2}\|\Lambda^{\beta}\tau^a\|_{L^2}\|\tau^a\|_{L^{\infty}}.
\end{align}
Applying $\Lambda$ to system (\ref{eq1}) and taking the $L^2$ inner product with $(\Lambda u^a,\Lambda\tau^a)$ yields
\begin{align*}
&\quad\frac{1}{2}\frac{d}{dt}\|(\Lambda u^a,\Lambda\tau^a)\|_{L^{2}}^{2}+\|\Lambda^{1+\beta}\tau^a\|_{L^{2}}^{2}+a\|\Lambda\tau^a\|_{L^{2}}^{2}\\
&=-\langle \Lambda (u^a\cdot\nabla u^a),\Lambda u^a\rangle
-\langle \Lambda (u^a\cdot\nabla\tau^a),\Lambda\tau^a\rangle
-\langle\Lambda Q(\nabla u^a,\tau^a),\Lambda\tau^a\rangle
\triangleq\sum_{i=1}^{3}I_i.
\end{align*}
Thanks to the two-dimensional cancellation, we have $I_1=0$. Applying Lemma \ref{lemma4} gives rise to
\begin{align*}
|I_2|&\lesssim \|\nabla u^a\|_{L^{\frac{2}{\beta}}}\|\nabla\tau^a\|_{L^2}
\|\nabla\tau^a\|_{L^{\frac{2}{1-\beta}}}\lesssim\|\Lambda^{2-\beta}u^a\|_{L^2}
\|\nabla\tau^a\|_{L^2}\|\Lambda^{1+\beta}\tau^a\|_{L^2}.
\end{align*}
By virtue of Lemma \ref{lemma5}, we have
\begin{align*}
|I_3|
&\lesssim\|\Lambda^{1+\beta}\tau^a\|_{L^2}
\left(\|\tau^a\|_{L^\infty}\|\Lambda^{2-\beta}u^a\|_{L^2}
+\|\Lambda^{1-\beta}\tau^a\|_{L^{\frac{2}{1-\beta}}}
 \|\nabla u^a\|_{L^{\frac2\beta}}\right)\\
&\lesssim\|\Lambda^{1+\beta}\tau^a\|_{L^2}\|\Lambda^{2-\beta}u^a\|_{L^2}
\left(\|\tau^a\|_{L^\infty}+\|\Lambda\tau^a\|_{L^2}\right),
\end{align*}
Putting together the estimates for $I_1$–$I_3$ and \eqref{211}, we obtain
\begin{align}\label{212}
&\quad\frac{1}{2}\frac{d}{dt}\|(u^a,\tau^a)\|_{H^{1}}^{2}
+\|\Lambda^{\beta}\tau^a\|_{H^{1}}^{2}+a\|\tau^a\|_{H^{1}}^{2}\lesssim \left(\|\tau^a\|_{H^1}+\|\tau^a\|_{L^\infty}\right)
\|\Lambda^{2-\beta} u^a\|_{L^{2}}\|\Lambda^\beta\tau^a\|_{H^1}.
\end{align}
To recover the dissipation of $u^a$, we consider the following cross-energy estimate. Let $k>0$ be a sufficiently small constant to be determined later. We then obtain
\begin{align}\label{214}
&\frac{d}{dt}\langle\Lambda^{1-\beta}u^a,k\Lambda^{1-\beta}\mathrm{div}\tau^a\rangle
+\frac{k}{2}\|\Lambda^{2-\beta}u^a\|_{L^2}^{2}\\
=&-k\langle\Lambda^{1-\beta}\mathbb P(u^a\cdot\nabla u^a),\Lambda^{1-\beta}\mathrm{div}\tau^a\rangle
-k\langle\Lambda^{1-\beta}\mathrm{div}(u^a\cdot\nabla\tau^a),\Lambda^{1-\beta}u^a\rangle\nonumber\\
&-k\langle\Lambda^{1-\beta}\mathrm{div}Q(\nabla u^a,\tau^a),\Lambda^{1-\beta}u^a\rangle
-k\langle\Lambda^{1-\beta}\Lambda^{2\beta}\mathrm{div}\tau^a,\Lambda^{1-\beta}u^a\rangle\nonumber\\\nonumber
&-ak\langle\Lambda^{1-\beta}\mathrm{div}\tau^a,\Lambda^{1-\beta}u^a\rangle
+k\|\Lambda^{1-\beta}\mathbb P\mathrm{div}\tau^a\|_{L^2}^{2}\nonumber\\\nonumber
\triangleq&\ \sum_{i=4}^{9}I_i,
\end{align}
where $\mathbb P$ denote the Leray projector.  Lemma \ref{lemma5} gives
\begin{align*}
|I_4|\lesssim k\|u^a\|_{L^{\frac{2}{2\beta-1}}}\|\Lambda^{3-3\beta}u^a\|_{L^{\frac{1}{1-\beta}}}\|\Lambda^{1+\beta}\tau^a\|_{L^2}
\lesssim k\|\Lambda^{2-2\beta}u^a\|_{L^2}\|\Lambda^{2-\beta}u^a\|_{L^2}\|\Lambda^{1+\beta}\tau^a\|_{L^2},\end{align*}
and
\begin{align*}
|I_5|&\lesssim k\|\Lambda^{2-\beta}u^a\|_{L^2}\Big(\|u^a\|_{L^{\frac{2}{2\beta-1}}}\|\Lambda^{1-\beta}\nabla\tau^a\|_{L^{\frac{1}{1-\beta}}}+\|\Lambda^{1-\beta}u^a\|_{L^{\frac{2}{1-\beta}}}\|\nabla\tau^a\|_{L^{\frac2\beta}}\Big)\\
&\lesssim k\|\Lambda^{2-\beta}u^a\|_{L^2}\Big(\|\Lambda^{2-2\beta}u^a\|_{L^{2}}\|\Lambda^{1+\beta}\tau^a\|_{L^{2}}+\|\Lambda u^a\|_{L^{2}}\|\Lambda^{2-\beta}\tau^a\|_{L^{2}}\Big)\\
&\lesssim k\|\Lambda^{2-\beta}u^a\|_{L^2}\| u^a\|_{H^1}\big(\|\Lambda^{1+\beta}\tau^a\|_{L^2}+\|\Lambda^{2-\beta}\tau^a\|_{L^2}\big),
\end{align*}
and
\begin{align*}
|I_6|&\lesssim k\|\Lambda^{2-\beta}u^a\|_{L^2}\Big(\|\tau^a\|_{L^\infty}\|\Lambda^{2-\beta}u^a\|_{L^2}+\|\Lambda^{1-\beta}\tau^a\|_{L^{\frac{2}{1-\beta}}}\|\nabla u^a\|_{L^{\frac2\beta}}\Big)\\*
&\lesssim k\Big(\|\tau^a\|_{L^\infty}+\|\Lambda\tau^a\|_{L^2}\Big)\|\Lambda^{2-\beta}u^a\|_{L^2}^2.
\end{align*}
The remaining three terms can be estimated as follows:
\begin{align*}
|I_7|+|I_8|+|I_9|&\lesssim k\|\Lambda^{1+\beta}\tau^a\|_{L^2}
\|\Lambda^{2-\beta}u^a\|_{L^2}+ak\|\Lambda^{1-\beta}\tau^a\|_{L^2}
\|\Lambda^{2-\beta}u^a\|_{L^2}+k\|\Lambda^{2-\beta}\tau^a\|_{L^2}^2.
\end{align*}
Young’s inequality then yields
\begin{align}\label{2130}
\sum_{i=4}^{9}|I_i|
\le&\ \frac{k}{4}\|\Lambda^{2-\beta}u^a\|_{L^2}^2
+Ck\big(\|u^a\|_{H^1}+\|\tau^a\|_{H^1}+\|\tau^a\|_{L^\infty}\big)\\\nonumber
&\quad\times\big(\|\Lambda^{2-\beta}u^a\|_{L^2}^2
+\|\Lambda^\beta\tau^a\|_{H^1}^2\big)
+Ck\|\Lambda^\beta\tau^a\|_{H^1}^2
+Cka^{2}\|\tau^a\|_{H^1}^2.
\end{align}
Combining \eqref{212}, \eqref{214} and \eqref{2130}, and choosing $k>0$ sufficiently small, we deduce that
\begin{align}\label{215}
&\frac12\frac{d}{dt}\left(\|(u^a,\tau^a)\|_{H^1}^2
+2k\langle\Lambda^{1-\beta}u^a,\Lambda^{1-\beta}\mathrm{div}\tau^a\rangle\right)+\frac{k}{8}\|\Lambda^{2-\beta}u^a\|_{L^2}^2
+\frac12\|\Lambda^\beta\tau^a\|_{H^1}^2
+\frac a2\|\tau^a\|_{H^1}^2\\\nonumber
\lesssim&\ \left(\|u^a\|_{H^1}+\|\tau^a\|_{H^1}
+\|\tau^a\|_{L^\infty}\right)
\Big(\|\Lambda^{2-\beta}u^a\|_{L^2}^2
+\|\Lambda^\beta\tau^a\|_{H^1}^2\Big).
\end{align}
Since $\frac{1}{2}<\beta<1$, we have
\begin{align*}
k|\langle\Lambda^{1-\beta}u^a,
\Lambda^{1-\beta}\mathrm{div}\tau^a\rangle|
&\lesssim k\|\Lambda^{2-2\beta}u^a\|_{L^2}
\|\Lambda\tau^a\|_{L^2}
\lesssim k\|(u^a,\tau^a)\|_{H^1}^2.
\end{align*}
In view of the bootstrap assumption \eqref{210}, inequality \eqref{215} yields, for any $0\le a\le1$,
\begin{align}\label{217}
&\frac12\frac{d}{dt}\left(\|(u^a,\tau^a)\|_{H^1}^2
+2k\langle\Lambda^{1-\beta}u^a,\Lambda^{1-\beta}\mathrm{div}\tau^a\rangle\right)\\\nonumber
&+\frac{k}{16}
\|\Lambda^{2-\beta}u^a\|_{L^2}^2
+\frac14\|\Lambda^\beta\tau^a\|_{H^1}^2
+\frac a2\|\tau^a\|_{H^1}^2\le 0.
\end{align}
Integrating \eqref{217} over $[0,t]$, we arrive at
\begin{align}\label{218}
\|(u^a(t),\tau^a(t))\|_{H^1}^2
\le C\|(u_{0},\tau_{0})\|_{H^1}^2.
\end{align}
We next derive a uniform $L^\infty$ bound for $\tau^a$. Duhamel's formula gives
\begin{align}\label{219}
\tau^a(t)={}&e^{-at}e^{-t(-\Delta)^\beta}\tau_0\\\nonumber
&-\int_0^t e^{-a(t-s)}e^{-(t-s)(-\Delta)^\beta}
(u^a\cdot\nabla\tau^a)(s)\,ds\\\nonumber
&+\int_0^t e^{-a(t-s)}e^{-(t-s)(-\Delta)^\beta}
\big(D(u^a)-Q(\nabla u^a,\tau^a)\big)(s)\,ds.\nonumber
\end{align}
Applying $\dot\Delta_j$ to \eqref{219} and using Lemma \ref{lem:heat-smoothing} together with Bernstein's inequality, we obtain, for $\frac1\beta<p<2$,
\begin{align}\label{220}
\|\sum_{j=0}^{\infty}\dot\Delta_{j}\tau^a(t)\|_{L^{\infty}}\lesssim\ &\|\tau_{0}\|_{L^{\infty}}
+\sum_{j=0}^{\infty}\int_{0}^{t}e^{-c(t-s)2^{2j\beta}}2^{\frac{2j}{p}}\|(u^a\cdot\nabla\tau^a)(s)\|_{L^{p}}ds\\\nonumber
&+\sum_{j=0}^{\infty}\int_0^te^{-c(t-s)2^{2j\beta}}\left(2^{\frac{2j}{p}}\|Q(\nabla u^a,\tau^a)(s)\big\|_{L^p}+2^{j}\|\nabla u^{a}(s)\|_{L^{2}}\right)ds\\\nonumber
\lesssim &\ \|\tau_{0}\|_{L^{\infty}}+\sum_{j=0}^{\infty}\int_{0}^{t}e^{-c(t-s)2^{2j\beta}}\left(2^{\frac{2j}{p}}\|(u^a,\tau^a)\|_{L^{\frac{2p}{2-p}}}\|\nabla(\tau^a,u^a)\|_{L^{2}}+2^{j}\|\nabla u^a\|_{L^2}\right)ds\\\nonumber
\lesssim &\ \|\tau_{0}\|_{L^{\infty}}+\sum_{j=0}^{\infty}\int_{0}^{t}e^{-c(t-s)2^{2j\beta}}\left(2^{\frac{2j}{p}}\|\Lambda^{2-\frac{2}{p}}(u^a,\tau^a)\|_{L^{2}}\|\nabla(\tau^a,u^a)\|_{L^{2}}+2^{j}\|\nabla u^a\|_{L^2}\right)ds\\\nonumber
\lesssim &\ \|\tau_{0}\|_{L^{\infty}}+\left(\|(u^a,\tau^a)\|_{L^{\infty}(0,t;H^{1})}^{2}+\|\nabla u^a\|_{L^{\infty}(0,t;L^{2})}\right)\sum_{j=0}^{\infty}\int_{0}^{t}e^{-c(t-s)2^{2j\beta}}2^{\frac{2j}{p}}ds.\nonumber
\end{align}
Hence,
\begin{align}\label{221}
\|\tau^a(t)\|_{L^{\infty}}&\le \|S_{0}\tau^{a}(t)\|_{L^{\infty}}+\|\sum_{j=0}^{\infty}\dot\Delta_{j}\tau^a(t)\|_{L^{\infty}}\\\nonumber
&\lesssim \|\tau^{a}(t)\|_{L^{2}}+\|\tau_{0}\|_{L^{\infty}}+\|(u^a,\tau^a)\|_{H^1}^{2}+\|\nabla u^a\|_{L^2}\\\nonumber
&\lesssim \|(u_{0},\tau_{0})\|_{H^{1}}+\|\tau_{0}\|_{L^{\infty}}+\|(u_0,\tau_0)\|_{H^1}^{2}.\nonumber
\end{align}
Together, \eqref{218} and \eqref{221} imply
\begin{align}\label{222}
\|(u^a(t),\tau^a(t))\|_{H^1}+\|\tau^a(t)\|_{L^{\infty}}
&\lesssim \|(u_{0},\tau_{0})\|_{H^{1}}+\|\tau_{0}\|_{L^{\infty}}+\|(u_0,\tau_0)\|_{H^1}^{2}.
\end{align}
The a priori bound now follows by a standard bootstrap argument.

\noindent\textbf{Proof of global existence of weak solutions for system (\ref{eq1}):}

Set $S_nf=\mathscr F^{-1}\left(\chi\left(\frac{\cdot}{n}\right)\mathscr Ff\right)$. We consider the following approximate system:
\begin{align}
\label{eq2}
\left\{
\begin{aligned}
&\partial_t u^n + S_n \mathbb{P}(v^n \cdot \nabla v^n) = \mathbb{P} \operatorname{div} \tau^n, \\
&\partial_t \tau^n + \Lambda^{2\beta} \tau^n + a \tau^n = D(u^n) - S_n(v^n \cdot \nabla \sigma^n + Q(\nabla v^n, \sigma^n)), \\
&(u^n, \tau^n)|_{t=0} = S_n(u_0, \tau_0),
\end{aligned}
\right.
\end{align}
where $v^n=S_nu^n$ and $\sigma^n=S_n\tau^n$. The Cauchy–Lipschitz theorem implies that (\ref{eq2}) admits a unique maximal solution $(u^n,\tau^n)$ in the space $C^1([0,T_n^*[;L^2\times L^2)$.
The above a priori estimates imply that $T_n^*=\infty$. For each $T>0$, the approximate solutions $(u^n,\tau^n)$ are uniformly bounded in
\begin{align}\label{223}
 L^\infty(0,T;H^1)\cap L^2(0,T;H^{2-\beta})\times
 L^\infty(0,T;H^1\cap L^\infty)
 \cap L^2(0,T;H^{1+\beta}).
\end{align}
By weak compactness, we obtain, up to a subsequence,
\begin{align*}
&u^n\overset{*}{\rightharpoonup}u
\quad\text{in }L^\infty(0,T;H^1),
\qquad
u^n\rightharpoonup u
\quad\text{in }L^2(0,T;H^{2-\beta}),\\
&\tau^n\overset{*}{\rightharpoonup}\tau
\quad\text{in }L^\infty(0,T;H^1),
\qquad
\tau^n\overset{*}{\rightharpoonup}\tau
\quad\text{in }L^\infty((0,T)\times\mathbb R^2),\\
&\tau^n\rightharpoonup\tau
\quad\text{in }L^2(0,T;H^{1+\beta}).
\end{align*}
It follows from (\ref{eq2}) that
\begin{align*}
\|\partial_tu^n\|_{H^{-1}}\lesssim \|u^n\|_{L^4}^2
 +\|\tau^n\|_{L^2}\lesssim \|u^n\|_{H^1}^2
 +\|\tau^n\|_{H^1}.
\end{align*}
Likewise, we have
\begin{align*}
\|\partial_t\tau^n\|_{H^{-1}}
&\lesssim \|v^n\sigma^n\|_{L^2}+\|\tau^n\|_{H^{2\beta-1}}+a\|\tau^n\|_{L^2}
 +\|\nabla u^n\|_{L^2}
   \big(1+\|\tau^n\|_{L^\infty}\big)\\
&\lesssim
 \|u^n\|_{H^1}\|\tau^n\|_{H^1}
 +\|\tau^n\|_{H^1}
 +\|u^n\|_{H^1}
   \big(1+\|\tau^n\|_{L^\infty}\big).
\end{align*}
Combining these estimates, we obtain
\begin{align}\label{224}
\sup_{n\ge1}
\big(\|\partial_tu^n\|_{L^\infty(0,T;H^{-1})}
+\|\partial_t\tau^n\|_{L^\infty(0,T;H^{-1})}\big)\le C_T.
\end{align}
Fix $T>0$ and let $K\subset\mathbb R^2$ be a bounded smooth domain. Since $\frac12<\beta<1$, the following embeddings hold:
\begin{align*}
H^{2-\beta}(K)\hookrightarrow\hookrightarrow H^1(K)\hookrightarrow H^{-1}(K),
\qquad
H^{1+\beta}(K)\hookrightarrow\hookrightarrow  H^1(K)\hookrightarrow H^{-1}(K).
\end{align*}
By the uniform estimates \eqref{223}–\eqref{224} and the Aubin–Lions lemma, we may extract a subsequence, still denoted by $(u^n,\tau^n)$, such that
\begin{align}\label{225}
(u^n,\tau^n)
\longrightarrow(u,\tau)
\quad\text{ in }
L^2(0,T;H^1(K)).
\end{align}
Taking $K=B_R$ with $R=1,2,\ldots$ and using a diagonal argument,
we obtain the convergence in
$L^2(0,T;H^1_{\mathrm{loc}}(\mathbb R^2))$.
Using the compact embedding
$H^1(K)\hookrightarrow\hookrightarrow H^\sigma(K)$ for any $0\le\sigma<1$, we also deduce
that
\begin{align*}
(u^n,\tau^n)
\longrightarrow(u,\tau)
\quad\text{in }C([0,T];H^\sigma(K)),\qquad 0\le\sigma<1.
\end{align*}
Since $\chi$ equals one in a neighborhood of the origin,
the uniform bounds in \eqref{223} imply
\begin{align*}
\|v^n-u^n\|_{L^2(0,T;H^1)}
&\lesssim n^{-(1-\beta)}
\|u^n\|_{L^2(0,T;H^{2-\beta})}\longrightarrow0,\\
\|\sigma^n-\tau^n\|_{L^2(0,T;H^1)}
&\lesssim n^{-\beta}
\|\tau^n\|_{L^2(0,T;H^{1+\beta})}\longrightarrow0.
\end{align*}
Combining these estimates with \eqref{225}, we obtain
\begin{align*}
(v^n,\sigma^n)\longrightarrow(u,\tau)
\quad\text{in }L^2(0,T;H^1_{\mathrm{loc}}(\mathbb R^2)).
\end{align*}
Consequently, for every $\phi\in  C_c^\infty([0,T)\times\mathbb R^2;\mathbb R^2)$with $\operatorname{div}\phi=0$, there holds
\begin{align*}
&\langle S_n(v^n\otimes v^n)-u\otimes u,\nabla\phi\rangle\\
&\quad=\langle v^n\otimes v^n-u\otimes u,\nabla\phi\rangle+v^n\otimes v^n,(S_n-I)\nabla\phi\rangle\longrightarrow0.
\end{align*}
The same argument applies to the $v^n\otimes\sigma^n$. Moreover, the bilinearity of $Q$ gives
\begin{align*}
&\|Q(\nabla v^n,\sigma^n)-Q(\nabla u,\tau)\|_{L^1((0,T)\times K)}
\\
&\quad\lesssim
\|\nabla v^n-\nabla u\|_{L^2((0,T)\times K)}
\|\sigma^n\|_{L^2((0,T)\times K)}
\\
&\qquad+
\|\nabla u\|_{L^2((0,T)\times K)}
\|\sigma^n-\tau\|_{L^2((0,T)\times K)}
\longrightarrow0.
\end{align*} Then, for every $\psi\in C_c^\infty((0,T)\times\mathbb R^2;\mathbb R^{2\times2})$, we conclude that
\begin{align*}
&\langle S_nQ(\nabla v^n,\sigma^n)-Q(\nabla u,\tau),\psi\rangle
\\
&\quad=
\langle Q(\nabla v^n,\sigma^n)-Q(\nabla u,\tau),\psi\rangle
+\langle Q(\nabla v^n,\sigma^n),(S_n-I)\psi\rangle
\longrightarrow0,
\end{align*}
Thus we can pass to the limit in all nonlinear terms of \eqref{eq2} in the sense of distributions.
\begin{align*}
&u\in L^\infty(0,\infty;H^1)
 \cap L^2(0,\infty; \dot H^{2-\beta}),\nonumber\\
&\tau\in L^\infty(0,\infty;H^1\cap L^\infty)
 \cap L^2(0,\infty;\dot H^{\beta}\cap\dot H^{1+\beta}).
\end{align*}
This establishes the existence of global weak solutions to system (\ref{eq1}).
\hfill$\Box$

\section{Optimal time decay rate with $a=0$}

	In this section, we consider the optimal time decay rate for (\ref{eq1}) with $a=0$. For the sake of convenience, we denote $u^0=u,\ \tau^0=\tau$ in this section. We shall present the proof of Theorem \ref{1theo2}. By a density argument, we only need to prove that the conclusions in this section hold for the smooth solutions.
    
We now turn to the large-time behavior of system \eqref{eq1}. Our analysis relies on the Fourier-splitting method of Schonbek \cite{Schonbek1985,Schonbek1991}.
To obtain additional dissipation for the velocity, we add a cross term to the basic $H^1$ energy functional, as in \cite{LLY}.
Motivated by the estimates established in Section 3, we define the energy and dissipation functionals as follows:
\begin{align*}
E_{0}={}&\|(u,\tau)\|_{H^{1}}^{2}
-2k\left\langle\Lambda^{2\beta-2}\nabla u,\tau
\right\rangle_{H^{2-2\beta}},\\
D_{0}={}&\frac{k}{16}\|\Lambda^{\beta}u\|_{H^{2-2\beta}}^{2}
+\frac18\|\Lambda^{\beta}\tau\|_{H^{1}}^{2},
\end{align*}
where $k>0$ will be chosen sufficiently small.
The corresponding energy inequality is stated below.
\begin{prop}{\label{3prop1}}
Let $a=0$, and let $(u,\tau)$ be a global weak solution obtained in Theorem \ref{1theo1}. Then there holds
\begin{align}\label{310}
\frac{d}{dt}E_0(t)+D_0(t)\le0.
\end{align}
\end{prop}
\begin{proof}
Arguing as in the derivation of (\ref{214}) yields
\begin{align}\label{311}
&\frac{d}{dt}\langle\Lambda^{\beta-1}u,k\Lambda^{\beta-1}\mathrm{div}\tau\rangle
+\frac{k}{2}\|\Lambda^\beta u\|_{L^2}^2\\\nonumber
=&-k\langle\Lambda^{\beta-1}\mathrm{div}(u\cdot\nabla\tau),\Lambda^{\beta-1}u\rangle-k\langle\Lambda^{\beta-1}\mathrm{div}Q(\nabla u,\tau),\Lambda^{\beta-1}u\rangle\\\nonumber
&
-k\langle\Lambda^{\beta-1}\Lambda^{2\beta}\mathrm{div}\tau,\Lambda^{\beta-1}u\rangle+k\|\Lambda^{\beta-1}\mathbb P\mathrm{div}\tau\|_{L^2}^2\\\nonumber
&-k\langle\Lambda^{\beta-1}\mathbb P(u\cdot\nabla u),\Lambda^{\beta-1}\mathrm{div}\tau\rangle\\\nonumber
\lesssim 
&\ k\|\Lambda^\beta u\|_{L^2}\left(
\|u\|_{L^{\frac2{1-\beta}}}\|\nabla\tau\|_{L^2}
+\|\nabla u\|_{L^2}\|\tau\|_{L^{\frac2{1-\beta}}}
+\|\Lambda^{3\beta-1}\tau\|_{L^2}\right)\\\nonumber
&+k\|\Lambda^\beta\tau\|_{L^2}^2+k\|u\|_{L^{\frac{2}{1-\beta}}}
\|\Lambda^\beta u\|_{L^{\frac2\beta}}
\|\Lambda^\beta\tau\|_{L^2}\\\nonumber
\lesssim&\ k\|\Lambda^\beta u\|_{L^2}\left(
\|\Lambda u\|_{L^2}\|\Lambda^\beta\tau\|_{L^2}
+\|\Lambda^\beta u\|_{L^2}\|\Lambda\tau\|_{L^2}
+\|\Lambda^{3\beta-1}\tau\|_{L^2}\right)+k\|\Lambda^\beta\tau\|_{L^2}^2\\\nonumber
\le&\ \frac{k}{4}\|\Lambda^\beta u\|_{L^2}^2
+Ck\|(u,\tau)\|_{H^1}\left(\|\Lambda^\beta u\|_{H^{2-2\beta}}^2
+\|\Lambda^\beta\tau\|_{H^1}^2\right)+Ck\|\Lambda^\beta\tau\|_{H^1}^2.
\end{align}
As $0<2\beta-1<1$, it follows that
\begin{align*}
\big|\langle\Lambda^{\beta-1}u,\Lambda^{\beta-1}\mathrm{div}\tau\rangle\big|
&\lesssim\|u\|_{L^2}\|\Lambda^{2\beta-1}\tau\|_{L^2}
\lesssim\|(u,\tau)\|_{H^1}^2.
\end{align*}
Adding \eqref{311} to \eqref{217} with $a=0$, and choosing $k$ and $\delta$ sufficiently small, we obtain \eqref{310}.
\end{proof}

Next we establish the $H^1$ decay estimate, which in particular yields the optimal zero-order $L^2$ decay.
\begin{prop}{\label{3prop2}}
  Under the same conditions as in Proposition \ref{3prop1}, if additionally $(u_0,\tau_0)\in \dot B_{2,\infty}^{-1}$, then there exists $C>0$ such that for every $t>0$, there holds
\begin{align*}
\|(u,\tau)(t)\|_{H^1}\lesssim(1+t)^{-\frac1{2\beta}},
\qquad t\ge0.
\end{align*}
\end{prop}
\begin{proof}
  Denote $S_{1}(t)=\left\{\xi:|\xi|^{2\beta}\le C_{2}(1+t)^{-1}\right\}$, $C_{2}$ is large enough. By Proposition \ref{3prop1}, we get
  \begin{align}\label{321}
  \frac{d}{dt}E_{0}(t)+\frac{C_{2}}{1+t}\left(k\|u\|_{H^1}^2+\|\tau\|_{H^1}^2\right)
  \lesssim \frac{1}{1+t}\int_{S_{1}(t)}
  \big(|\widehat{u}(\xi)|^{2}+|\widehat\tau(\xi)|^{2}\big)\,d\xi.
  \end{align}
  For convenience, we set
\[
F=-u\cdot\nabla u,\qquad
G=-u\cdot\nabla\tau-Q(\nabla u,\tau).
\]
Applying the Fourier transform to \eqref{eq1}, we have
\begin{align}\label{322}
\left\{\begin{array}{l}\widehat{u}^{j}_{t}+i\xi_{j}\widehat{P}-i\xi_{k}\widehat{\tau}^{ jk}=\widehat{F}^{j},\\ [1ex] \widehat{\tau}^{jk}_{t}+|\xi|^{2\beta}\widehat{\tau}^{jk}-\frac{i}{2}(\xi_{k}\widehat{u}^{j}+ \xi_{j}\widehat{u}^{k})=\widehat{G}^{jk}.
\end{array}\right.
\end{align}
Taking the inner product with $(\widehat u,\widehat\tau)$, we obtain
\begin{align}\label{323}
\frac{1}{2}\frac{d}{dt}(| \widehat{u}|^{2}+|\widehat{\tau}|^{2})+|\xi|^{2\beta}|\widehat{\tau}|^{2}=\mathcal{R}e[\widehat{F} \cdot\overline{\widehat{u}}]+\mathcal{R}e[\widehat{G}:\overline{\widehat{\tau}}].
\end{align}
Integrating (\ref{323}) in time on $[0,t]$, we get
\begin{align}\label{324}
|\widehat{u}|^{2}+|\widehat{\tau}|^{2}\lesssim (|\widehat{u_{0}}|^{2}+|\widehat{\tau_0}|^{2})+\int_{0}^{t}|\widehat{F}\cdot\overline{\widehat{u}}|+|\widehat{G}:\overline{\widehat{\tau}}|ds.
\end{align}
Integrating (\ref{324}) over $S_{1}(t)$ with $\xi$, then we derive
\begin{align}\label{325}
\int_{S_{1}(t)}|\widehat{u}|^{2}+|\widehat{\tau}|^{2}d\xi\lesssim \int_{S_{1}(t)}|\widehat{u_{0}}|^{2}+|\widehat{\tau_{0}}|^{2}d\xi+\int_{S_{1}(t)}\int_{0}^{t}|\widehat{F}\cdot\overline{\widehat{u}}|+|\widehat{G}:\overline{\widehat{\tau}}|\,ds\,d\xi.
\end{align}
Since the initial energy is finite and
$(u_0,\tau_0)\in\dot B^{-1}_{2,\infty}$, we have
\begin{align}\label{326}
\int_{S_{1}(t)}(|\widehat{u_{0}}|^{2}+|\widehat{\tau_{0}}|^{2})d\xi &\lesssim\sum_{j\leq\log_{2}\left[\frac{4}{3}C_{2}^{\frac{1}{2\beta}}(1+t)^{-\frac{1}{2\beta}}\right]}\int_{\mathbb{R}^{2}}2\varphi^{2}(2^{-j}\xi)\left(|\widehat{u_{0}}|^{2}+|\widehat{\tau_{0}}|^{2}\right)d\xi \\ \nonumber
 & \lesssim\sum_{j\leq\log_{2}\left[\frac{4}{3}C_{2}^{\frac{1}{2\beta}}(1+t)^{-\frac{1}{2\beta}}\right]}\left(\|\dot{\Delta}_{j}u_{0}\|_{L^{2}}^{2}+\|\dot{\Delta}_{j}\tau_{0}\|_{L^{2}}^{2}\right) \\ \nonumber
 & \lesssim (1+t)^{-\frac{1}{\beta}}\|(u_{0},\tau_{0})\|_{\dot{B}_{2,\infty}^{-1}}^{2}.\nonumber
\end{align}
Minkowski's inequality gives
\begin{align}\label{327}
\int_{S_{1}(t)}\int_{0}^{t}|\widehat{F}\cdot\overline{\widehat{u}}|+|\widehat{G}:\overline{\widehat{\tau}}|\,ds\,d\xi=& \int_{0}^{t}\int_{S_{1}(t)}|\widehat{F}\cdot\overline{\widehat{u}}|+|\widehat{G}:\overline{\widehat{\tau}}|d\xi ds \\\nonumber
\lesssim &\ |S_{1}(t)|^{\frac{1}{2}}\int_{0}^{t}\|F\|_{L^1}\|u\|_{L^2}+\|G\|_{L^1}\|\tau\|_{L^2}ds\\\nonumber
\lesssim &\ (1+t)^{-\frac{1}{2\beta}}\int_{0}^{t}\left(\|u\|_{L^2}^2+\|\tau\|_{L^2}^2\right)\left(\|\nabla u\|_{L^2}+\|\nabla\tau\|_{L^2}\right)ds\\\nonumber
\lesssim &\ (1+t)^{-\frac{1}{2\beta}+\frac{1}{2}}\left(\int_{0}^{t}D_{0}(s)ds\right)^{\frac{1}{2}}\\\nonumber
\lesssim &\ (1+t)^{-\frac{1}{2\beta}+\frac{1}{2}}.\nonumber
\end{align}
Substituting (\ref{326}) and (\ref{327}) into (\ref{325}), we have
\begin{align}\label{328}
\int_{S_1(t)}
\left(|\widehat{u}(\xi)|^2+
|\widehat{\tau}(\xi)|^2\right)\,d\xi
\lesssim (1+t)^{-\frac{1}{2\beta}+\frac12}.
\end{align}
Plugging (\ref{328}) into (\ref{321}) gives rise to
\begin{align*}
  \frac{d}{dt}E_{0}(t)+\frac{C_{2}}{1+t}\left(k\|u\|_{H^1}^2+\|\tau\|_{H^1}^2\right)
  \lesssim (1+t)^{-\frac1{2\beta}-\frac12}.
  \end{align*}
This establishes the preliminary decay bound
\begin{align}\label{3210}
E_{0}(t)\lesssim (1+t)^{-\frac{1}{2\beta}+\frac{1}{2}}.
  \end{align}

To improve the preliminary decay bound in \eqref{3210}, we use a time-weighted energy estimate. Recalling (\ref{324})--(\ref{327}) and using the equivalence of
$E_0$ and $\|(u,\tau)\|_{H^1}^2$, we infer from (\ref{321}) that
\begin{align}\label{3211}
\frac{d}{dt}E_0(t)+\frac{C_2}{1+t}\left(k\|u\|_{H^1}^2+\|\tau\|_{H^1}^2\right)
&\lesssim (1+t)^{-1-\frac1\beta}+(1+t)^{-1-\frac1{2\beta}}
  \int_0^t E_0(s)D_0(s)^{\frac12}\,ds.
\end{align}
Multiplying (\ref{3211}) by $(1+t)^{1+\frac{1}{2\beta}}$ and then integrating over $[0,t]$. we obtain
\begin{align}\label{3212}
E_0(t)&\lesssim (1+t)^{-\frac1\beta}
 +(1+t)^{-\frac1{2\beta}}
 \int_0^t E_0(s)D_0(s)^{\frac12}\,ds.
\end{align}
Set
$N(t)=\sup_{0\le s\le t}(1+s)^{\frac1{2\beta}}E_0(s)$.
By (\ref{3212}), we have
\begin{align*}
N(t)&\lesssim 1+\int_0^t (1+s)^{-\frac1{2\beta}}
 N(s)D_0(s)^{\frac12}\ ds.
\end{align*}
An application of Gronwall’s inequality yields
\begin{align*}
N(t)&\lesssim \mathrm{exp}\left(C\displaystyle\int_0^\infty(1+s)^{-\frac1{2\beta}}D_0(s)^{\frac12}ds\right)\\
&\lesssim \mathrm{exp}\left(\left(C\int_0^\infty(1+s)^{-\frac1\beta}\,ds\right)^{\frac12}
\left(\int_0^\infty D_0(s)\,ds\right)^{\frac12}\right)<\infty.
\end{align*}
In particular,
\begin{align}\label{3222}
E_0(t)\lesssim (1+t)^{-\frac1{2\beta}}.
\end{align}
Fix $1-\frac1{2\beta}<\delta<\frac1{2\beta}$.
Multiplying (\ref{310}) by $(1+t)^\delta$, integrating on
$[0,t]$, and using (\ref{3222}), we obtain
\begin{align}\label{3223}
&(1+t)^\delta E_0(t)
+\int_0^t(1+s)^\delta D_0(s)\,ds\\
&\qquad\lesssim E_0(0)
+\int_0^t(1+s)^{\delta-1}E_0(s)\,ds\nonumber\\
&\qquad\lesssim E_0(0)
+\int_0^\infty(1+s)^{\delta-1-\frac1{2\beta}}\,ds
\lesssim1.\nonumber
\end{align}
We next show that the regularity in
$\dot B^{-1}_{2,\infty}$ is propagated for all time. Applying $\dot\Delta_{j}$ to (\ref{eq1}), we find
\begin{align*}
 \left\{\begin{array}{l}\dot{\Delta}_{{j}}u_{t}+ \nabla\dot{\Delta}_{j}P-\text{div}\ \dot{\Delta}_{j}\tau=\dot{\Delta}_{j}F,\\ [1ex] \dot{\Delta}_{j}\tau_{t}+(-\Delta)^{\beta}\dot{\Delta}_{j}\tau-\dot{\Delta}_{j}D(u)=\dot{\Delta}_{j}G.\end{array}\right.
\end{align*}
Then we get
\begin{align}\label{3226}
\frac{d}{dt}(\|\dot{\Delta}_{j}u\|_{L^{2}}^{2}+\|\dot{\Delta}_{j}\tau\|_{L^{2} }^{2})+2\|\Lambda^{\beta}\dot{\Delta}_{j}\tau\|_{L^{2}}^{2}\lesssim \|\dot{\Delta}_{j}F\| _{L^{2}}\|\dot{\Delta}_{j}u\|_{L^{2}}+\|\dot{\Delta}_{j}G\|_{L^{2}}\|\dot{ \Delta}_{j}\tau\|_{L^{2}}.
\end{align}
Multiplying (\ref{3226}) by $2^{-2j}$ and taking $l^{\infty}$ norm, we have
\begin{align*}
\frac{d}{dt}\left(\|u\|_{\dot{B}^{-1}_{2,\infty}}^{2}+\|\tau\|_{\dot{B}^{-1}_{2,\infty}}^{2}\right)\lesssim \|F\|_{\dot{B}^{-1}_{2,\infty}}\|u\|_{ \dot{B}^{-1}_{2,\infty}}+\|G\|_{\dot{B}^{-1}_{2,\infty}}\|\tau\|_{ \dot{B}^{-1}_{2,\infty}}.
\end{align*}
Let $\displaystyle M(t)=\sup_{0\le s\le t}\left(\|u\|_{\dot{B}^{-1}_{2,\infty}}+\|\tau\|_{\dot{B}^{-1}_{2,\infty}}\right)$, then we yield
\begin{align*}
M^{2}(t)\lesssim M^{2}(0)+M(t)\int_{0}^{t}\left(\|F\|_{\dot{B}^{-1}_{2,\infty}}+\|G\|_{\dot{B}^{-1}_{2,\infty}}\right)ds.
\end{align*}
By $L^1\hookrightarrow\dot B^{-1}_{2,\infty}$ and
(\ref{3223}), we deduce
\begin{align}\label{3229}
\int_0^t\big(\|F(s)\|_{\dot B^{-1}_{2,\infty}}
+\|G(s)\|_{\dot B^{-1}_{2,\infty}}\big)\,ds
&\lesssim\int_0^t(\|F(s)\|_{L^1}+\|G(s)\|_{L^1})\,ds\\
&\lesssim\int_0^t\|(u,\tau)(s)\|_{L^2}
 \|\nabla(u,\tau)(s)\|_{L^2}\,ds\nonumber\\
&\lesssim1.\nonumber
\end{align}
Thus, $M(t)\le C$. We now estimate the nonlinear
contribution using this bound together with
\eqref{3222} and \eqref{3229}:
\begin{align}\label{3230}
\int_{S_{1}(t)}\int_{0}^{t}|\widehat{F}\cdot\overline{\widehat{u}}|+|\widehat{G}:\overline{\widehat{\tau}}|\,ds\,d\xi=& \int_{0}^{t}\int_{S_{1}(t)}|\widehat{F}\cdot\overline{\widehat{u}}|+|\widehat{G}:\overline{\widehat{\tau}}|d\xi ds \\\nonumber
\lesssim & \int_{0}^{t}\left(\|F\|_{L^1}\int_{S_{1}(t)}|{\widehat{u}}|d\xi+\|G\|_{L^1}\int_{S_{1}(t)}|{\widehat{\tau}}|d\xi\right
) ds \\\nonumber
\lesssim & \ (1+t)^{-\frac{1}{2\beta}}\int_{0}^{t}\left[\left(\|F\|_{L^1}+\|G\|_{L^1}\right)\left(\int_{S_{1}(t)}|\widehat{u}|^{2}+|\widehat{\tau}|^{2}d\xi\right)^{\frac{1}{2}}\right]ds\\\nonumber
\lesssim &\ (1+t)^{-\frac1\beta}M(t)
\int_0^t\big(\|F(s)\|_1+\|G(s)\|_1\big)\,ds\\\nonumber
\lesssim &\ (1+t)^{-\frac1\beta}.\nonumber
\end{align}
It follows from \eqref{321}, \eqref{326}, and
\eqref{3230} that
\begin{align}\label{3231}
E_0(t)\lesssim (1+t)^{-\frac1\beta}.
\end{align}
 \end{proof}

\begin{prop}\label{3prop3}
Under the same conditions as in Proposition \ref{3prop2}, then there exists $C>0$ such that for every $t>0$, there holds
\begin{align}\label{3234}
\|\tau(t)\|_{L^\infty}\le C(1+t)^{-\frac1{2\beta}}.
\end{align}
\end{prop}
\begin{proof}
We write
\begin{align*}
\tau(t)=S_{0}\tau(t)+\sum_{j\ge0}\dot\Delta_{j}\tau(t).
\end{align*}
Using (\ref{3231}), we deduce that
\begin{align}\label{3235}
\|S_{0}\tau(t)\|_{L^\infty}\lesssim\|S_{0}\tau(t)\|_{L^2}
\le C(1+t)^{-\frac1{2\beta}}.
\end{align}
 For the high-frequency part, the Duhamel formula takes the form
\begin{align*}
\sum_{j\ge0}\dot\Delta_j\tau(t)={}&\sum_{j\ge0}e^{-t(-\Delta)^\beta}\dot\Delta_j\tau_0\\
&+\sum_{j\ge0}\int_0^t e^{-(t-s)(-\Delta)^\beta}
\dot\Delta_j\left(D(u)-u\cdot\nabla\tau-Q(\nabla u,\tau)\right)(s)\,ds.\nonumber
\end{align*}
By Lemma \ref{lem:heat-smoothing} and Bernstein's inequality, the initial data term is bounded by
\begin{align}\label{3237}
\|\sum_{j\ge0}e^{-t(-\Delta)^{\beta}}\dot\Delta_{j}\tau_0\|_{L^\infty}&\lesssim \sum_{j\in\mathbb Z}2^je^{-ct2^{2\beta j}}\|\dot\Delta_j\tau_0\|_{L^2}\\
&\lesssim \|\tau_0\|_{\dot B^{-1}_{2,\infty}}\sum_{j\in\mathbb Z}2^{2j}e^{-ct2^{2\beta j}}\nonumber\\\nonumber
&\lesssim t^{-\frac1\beta}\|\tau_0\|_{\dot B^{-1}_{2,\infty}},
\end{align}
Repeating the calculation in the proof of \eqref{220} with the estimates \eqref{3231} and \eqref{3237}, we infer, for $t>1$, 
\begin{align}\label{3238}
\|\sum_{j\ge0}\dot\Delta_{j}\tau(t)\|_{L^{\infty}}\lesssim\ & t^{-\frac{1}{\beta}}+\sum_{j\ge0}\int_0^t
e^{-c(t-s)2^{2\beta j}}
\left(
2^{\frac{2j}{p}}(1+s)^{-\frac{1}{\beta}}
+2^j(1+s)^{-\frac{1}{2\beta}}
\right)\,ds\\\nonumber
\lesssim{}&
(1+t)^{-\frac{1}{2\beta}}
+\sum_{j\ge0}
\left(2^{\frac{2j}{p}}+2^j\right)
\int_0^{t/2}
e^{-c(t-s)2^{2\beta j}}\,ds\\\nonumber
&+(1+t)^{-\frac{1}{2\beta}}
\sum_{j\ge0}
\left(2^{\frac{2j}{p}}+2^j\right)
\int_{t/2}^t
e^{-c(t-s)2^{2\beta j}}\,ds\\\nonumber
\lesssim{}&
(1+t)^{-\frac{1}{2\beta}}
+\sum_{j\ge0}
\left(2^{\frac{2j}{p}}+2^j\right)
2^{-2\beta j}e^{-ct2^{2\beta j}/2}\\\nonumber
&+(1+t)^{-\frac{1}{2\beta}}
\sum_{j\ge0}
\left(2^{\frac{2j}{p}}+2^j\right)2^{-2\beta j}\\
\lesssim{}&
(1+t)^{-\frac{1}{2\beta}}.\nonumber
\end{align}
Combining (\ref{3235}) and (\ref{3238}), we get for $t\ge 1$,
\begin{align*}
\|\tau(t)\|_{L^\infty}\le C(1+t)^{-\frac1{2\beta}}.
\end{align*}
For $0\le t\le 1$, the same estimate follows from the uniform $L^\infty$ bound in (\ref{222}). Consequently,
\[
\|\tau(t)\|_{L^\infty}
\lesssim (1+t)^{-\frac{1}{2\beta}},
\qquad t\ge0.
\]
This proves (\ref{3234}).
\end{proof}

To establish the optimal decay rate for the first-order derivatives, we introduce the following weighted energy and dissipation functionals:
 \begin{align*}
 \widetilde{E}_{1}&=(1+t)^{a^{\prime}}\|\Lambda(u,\tau)\|_{L^2}^{2}
 +k\langle\Lambda^{1-\beta}\tau,-\nabla\Lambda^{1-\beta}u\rangle,\\
 \widetilde{D}_{1}&=(1+t)^{a^{\prime}}\|\Lambda^{1+\beta}\tau\|_{L^2}^{2}
 +\frac{k}{4}\|\Lambda^{2-\beta}u\|_{L^2}^{2},
 \end{align*}
 where $a^{\prime}=2-\frac{1}{\beta}\in(0,1)$ and $k>0$ is sufficiently small.
 Next we prove the optimal $\dot H^1$ time decay rate.
 \begin{prop}\label{3prop4}
 Under the same conditions as in Proposition \ref{3prop3}, there exists $C>0$ such that, for every $t\ge0$,
 \begin{align*}
 \|\Lambda^{s_1}(u,\tau)(t)\|_{L^2}\le C(1+t)^{-\frac{1+s_1}{2\beta}},
 \qquad 0\le s_1\le1.
 \end{align*}
 \end{prop}
 \begin{proof}
 A direct computation using \eqref{eq1} shows that
 \begin{align}\label{3335}
 \frac{d}{dt}\widetilde{E}_{1}+2\widetilde{D}_{1}
 &=a^{\prime}(1+t)^{a^{\prime}-1}\|\Lambda(u,\tau)\|_{L^2}^{2}
 +2(1+t)^{a^{\prime}}\langle\Lambda G,\Lambda\tau\rangle\\
 &\quad+k\langle-\Lambda^{1-\beta}G+\Lambda^{1+\beta}\tau,\nabla\Lambda^{1-\beta}u\rangle
 +k\langle\Lambda^{1-\beta}\mathbb{P}(F+\operatorname{div}\tau),\Lambda^{1-\beta}\operatorname{div}\tau\rangle.\nonumber
 \end{align}
 Using (\ref{212}), (\ref{3231}) and (\ref{3234}), we obtain
 \begin{align}\label{3337}
 (1+t)^{a^{\prime}}|\langle\Lambda G,\Lambda\tau\rangle|\lesssim &\ (1+t)^{a^{\prime}}\big(\|\nabla\tau\|_{L^2}+\|\tau\|_{L^\infty}\big)
 \|\Lambda^{2-\beta}u\|_{L^2}\|\Lambda^{1+\beta}\tau\|_{L^2}\\
 \lesssim&\  \delta^{1-\beta}(1+t)^{a^{\prime}-\frac12}\|\Lambda^{2-\beta}u\|_{L^2}\|\Lambda^{1+\beta}\tau\|_{L^2}\nonumber\\
\le&\ \frac18(1+t)^{a^{\prime}}\|\Lambda^{1+\beta}\tau\|_{L^2}^{2}+C\delta^{2(1-\beta)}(1+t)^{a^{\prime}-1}\|\Lambda^{2-\beta}u\|_{L^2}^{2}\nonumber\\
\le&\ \frac18(1+t)^{a^{\prime}}\|\Lambda^{1+\beta}\tau\|_{L^2}^{2}
 +\frac{k}{32}\|\Lambda^{2-\beta}u\|_{L^2}^{2}.\nonumber
 \end{align}
 Similarly, we have
 \begin{align}\label{3338}
 &k\big|\langle\Lambda^{1-\beta}\mathbb{P}F,\Lambda^{1-\beta}\operatorname{div}\tau\rangle\big|
 +k\big|\langle\Lambda^{1-\beta}\operatorname{div}(u\cdot\nabla\tau),\Lambda^{1-\beta}u\rangle\big|\\\nonumber
 \lesssim&\ k\|\Lambda^{2-\beta}u\|_{L^2}
 \Big(\|\Lambda^{2-2\beta}u\|_{L^2}\|\Lambda^{1+\beta}\tau\|_{L^2}
 +\|\Lambda u\|_{L^2}\|\Lambda^{2-\beta}\tau\|_{L^2}\Big)\\\nonumber
 \le&\ \frac{k}{32}\|\Lambda^{2-\beta}u\|_{L^2}^{2}+Ck\delta^2\big(\|\Lambda^{1+\beta}\tau\|_{L^2}^{2} +\|\Lambda^{2-\beta}\tau\|_{L^2}^{2}\big),\nonumber
 \end{align}
 and
 \begin{align}\label{3339}
 k\big|\langle\Lambda^{1-\beta}\operatorname{div}Q(\nabla u,\tau),\Lambda^{1-\beta}u\rangle\big|
 \lesssim k\big(\|\tau\|_{L^\infty}+\|\Lambda\tau\|_{L^2}\big)\|\Lambda^{2-\beta}u\|_{L^2}^{2}\le\frac{k^{2}}{64}\|\Lambda^{2-\beta}u\|_{L^2}^{2},
 \end{align}
 and
 \begin{align}\label{3340}
&k\big|\langle\Lambda^{1+\beta}\tau,\nabla\Lambda^{1-\beta}u\rangle\big|+k\|\Lambda^{1-\beta}\mathbb P\operatorname{div}\tau\|_{L^2}^{2}\\
 &\le\frac{k}{64}\|\Lambda^{2-\beta}u\|_{L^2}^{2}+Ck\|\Lambda^{1+\beta}\tau\|_{L^2}^{2}+Ck\|\Lambda^{2-\beta}\tau\|_{L^2}^{2}.\nonumber
 \end{align}
 The interpolation inequality and \eqref{3231} imply
 \begin{align}\label{3341}
 \|\Lambda^{2-\beta}\tau\|_{L^2}^{2}
 &\lesssim\|\tau\|_{L^2}^{\frac{2(2\beta-1)}{1+\beta}}
 \|\Lambda^{1+\beta}\tau\|_{L^2}^{\frac{2(2-\beta)}{1+\beta}}\\
 &\le\frac{1}{100}(1+t)^{a^{\prime}}\|\Lambda^{1+\beta}\tau\|_{L^2}^{2}
 +C(1+t)^{-\frac{a^{\prime}(2-\beta)}{2\beta-1}}\|\tau\|_{L^2}^{2}\nonumber\\
&\lesssim\frac{1}{100}(1+t)^{a^{\prime}}\|\Lambda^{1+\beta}\tau\|_{L^2}^{2}
 +C(1+t)^{1-\frac3\beta}.\nonumber
 \end{align}
Plugging from (\ref{3337}) to (\ref{3340}) into (\ref{3335}), using (\ref{3341}) and taking $k$ sufficiently small, we deduce
 \begin{align*}
 \frac{d}{dt}\widetilde E_1
 +\frac12(1+t)^{a^{\prime}}\|\Lambda^{1+\beta}\tau\|_{L^2}^{2}
 +\frac{k}{4}\|\Lambda^{2-\beta}u\|_{L^2}^{2}
 \le a^{\prime}(1+t)^{a^{\prime}-1}\|\Lambda(u,\tau)\|_{L^2}^{2}
 +C(1+t)^{1-\frac3\beta}.
 \end{align*}
 Recall that $S_1(t)=\left\{\xi:|\xi|^{2\beta}\le C_2(1+t)^{-1}\right\}$. Then we obtain
 \begin{align*}
 \frac{d}{dt}\widetilde E_1
 +\Big(\frac{kC_2}{4}-a^{\prime}\Big)(1+t)^{a^{\prime}-1}\|\Lambda(u,\tau)\|_{L^2}^{2}
 \lesssim(1+t)^{1-\frac3\beta},
 \end{align*}
 where $\frac{kC_2}{4}-a^{\prime}>\frac3\beta-1$.
 Multiplying this inequality by $(1+t)^{\frac{3}{\beta}-1}$
and integrating over $[0,t]$, we deduce that
 \begin{align*}
 &(1+t)^{\frac3\beta-1}\widetilde E_1(t)
 +\Big(\frac{kC_2}{4}-a^{\prime}-\frac3\beta+1\Big)\int_0^t(1+s)^{\frac2\beta}\|\Lambda(u,\tau)(s)\|_{L^2}^{2}\,ds\\
 &\quad\lesssim(1+t)
 +k\Big(\frac3\beta-1\Big)\int_0^t(1+s)^{\frac3\beta-2}
 \big|\langle\Lambda^{1-\beta}\tau,-\nabla\Lambda^{1-\beta}u\rangle\big|\,ds.
 \end{align*}
 For any $\varepsilon>0$, the remaining cross term can be estimated as
 \begin{align*}
 (1+t)^{\frac3\beta-2}\big|\langle\Lambda^{1-\beta}\tau,-\nabla\Lambda^{1-\beta}u\rangle\big|\lesssim&\ (1+t)^{\frac3\beta-2}
 \|\Lambda(u,\tau)\|_{L^2}^{3-2\beta}\|(u,\tau)\|_{L^2}^{2\beta-1}\\
\le&\ \varepsilon(1+t)^{\frac2\beta}\|\Lambda(u,\tau)\|_{L^2}^{2}+C_\varepsilon\|(u,\tau)\|_{L^2}^{2}\\
\lesssim& \ \varepsilon(1+t)^{\frac2\beta}\|\Lambda(u,\tau)\|_{L^2}^{2}
 +C_\varepsilon(1+t)^{-\frac1\beta},
 \end{align*}
which implies that
 \begin{align*}
 (1+t)^{\frac3\beta-1}\widetilde E_1(t)
 +c\int_0^t(1+s)^{\frac2\beta}\|\Lambda(u,\tau)(s)\|_{L^2}^{2}\,ds
 \lesssim 1+t.
 \end{align*}
 Applying the same interpolation inequality once again, we find
 \begin{align*}
 (1+t)^{1+\frac2\beta}\|\Lambda(u,\tau)(t)\|_{L^2}^{2}&\lesssim(1+t)^{\frac3\beta-1}\widetilde E_1(t)
 +k(1+t)^{\frac3\beta-1}\big|\langle\Lambda^{1-\beta}\tau,-\nabla\Lambda^{1-\beta}u\rangle\big|\\
 &\le(1+t)^{\frac3\beta-1}\widetilde E_1(t)
 +\frac12(1+t)^{1+\frac2\beta}\|\Lambda(u,\tau)(t)\|_{L^2}^{2}
 +C(1+t)\|(u,\tau)\|_{L^2}^{2}\\
&\lesssim(1+t)+\frac12(1+t)^{1+\frac2\beta}\|\Lambda(u,\tau)\|_{L^2}^{2}.
 \end{align*}
 Consequently,
 \begin{align}\label{33410}
 \|\Lambda(u,\tau)(t)\|_{L^2}\lesssim(1+t)^{-\frac1\beta}.
 \end{align}
 Interpolating between \eqref{3231} and \eqref{33410} completes the proof for $0\le s_1\le1$.
 \end{proof}

\begin{prop}\label{3prop5}
 Under the same conditions as in Proposition \ref{3prop3}, there exists $C>0$ such that, for every $t\ge0$,
 \begin{align*}
\|\tau\|_{L^{\infty}}\le C(1+t)^{-\frac{1}{\beta}}.
 \end{align*}
\end{prop}
\begin{proof}
We first consider $t\ge1$. We write
\[
\tau=\Delta_{-1}\tau+\sum_{j\ge0}\dot\Delta_j\tau.
\]
By Proposition \ref{3prop4}, we have
\begin{equation}\label{3342}
\begin{aligned}
\|(F,G)\|_{L^1}
&\lesssim\|(u,\tau)\|_{L^2}\|\nabla(u,\tau)\|_{L^2}
\lesssim(1+s)^{-\frac{3}{2\beta}}.
\end{aligned}
\end{equation}
We begin with the low frequencies. Applying $\Delta_{-1}$ to (\ref{322}), we get
\begin{align*}
&\frac12\frac{d}{dt}\left(
|\widehat{\Delta_{-1}u}|^2+|\widehat{\Delta_{-1}\tau}|^2\right)
+|\xi|^{2\beta}|\widehat{\Delta_{-1}\tau}|^2=\operatorname{Re}\left(
\widehat{\Delta_{-1}F}\cdot\overline{\widehat{\Delta_{-1}u}}
+\widehat{\Delta_{-1}G}:\overline{\widehat{\Delta_{-1}\tau}}\right).
\end{align*}
The same argument as in the derivation of
\eqref{311} leads to
\[
\begin{aligned}
&\frac{d}{dt}\left[|\xi|^{2\beta-2}\operatorname{Re}\left(
\widehat{\Delta_{-1}u}\cdot
\overline{\widehat{\Delta_{-1}\operatorname{div}\tau}}\right)\right]
+\frac14|\xi|^{2\beta}|\widehat{\Delta_{-1}u}|^2\\
&\qquad\lesssim |\xi|^{2\beta}|\widehat{\Delta_{-1}\tau}|^2
+\left(|\widehat{\Delta_{-1}u}|+|\widehat{\Delta_{-1}\tau}|\right)
\left(|\widehat{\Delta_{-1}F}|+|\widehat{\Delta_{-1}G}|\right).
\end{aligned}
\]
Here we used Young's inequality and the fact that $\mathrm{supp}\widehat{\Delta_{-1}f}
\subset B(0,4/3)$. Let
\[
 E_{-1}(t)=|\widehat{\Delta_{-1}u}|^2
+|\widehat{\Delta_{-1}\tau}|^2
+2k|\xi|^{2\beta-2}\operatorname{Re}
\left(\widehat{\Delta_{-1}u}\cdot
\overline{\widehat{\Delta_{-1}\operatorname{div}\tau}}\right).
\]
Taking $k>0$ sufficiently small, then
we arrive at
\begin{align*}
\frac{d}{dt}E_{-1}+c|\xi|^{2\beta}E_{-1}
\lesssim E_{-1}^{\frac12}
\left(|\widehat{\Delta_{-1}F}|+|\widehat{\Delta_{-1}G}|\right),
\end{align*}
 which implies that
\begin{equation}\label{3344}
\begin{aligned}
|\widehat{\Delta_{-1}\tau}(t,\xi)|
\lesssim{}&e^{-ct|\xi|^{2\beta}}
|\widehat{\Delta_{-1}(u_0,\tau_0)}(\xi)|\\
&+\int_0^t e^{-c(t-s)|\xi|^{2\beta}}
\left(|\widehat{\Delta_{-1}F}(s,\xi)|
+|\widehat{\Delta_{-1}G}(s,\xi)|\right)\,ds.
\end{aligned}
\end{equation}
As in the proof of \eqref{326}, the contribution
from the initial data satisfies
\[
\int_{|\xi|\le \frac{4}{3}}
e^{-ct|\xi|^{2\beta}}\left(|\widehat u_0|+|\widehat\tau_0|\right)\,d\xi
\lesssim (1+t)^{-\frac{1}{\beta}}\|(u_0,\tau_0)\|_{\dot B^{-1}_{2,\infty}}.
\]
By virtue of \eqref{3342}, \eqref{3344} and the assumption
$\frac12<\beta<1$, the low-frequency part satisfies
\begin{equation}\label{3345}
\begin{aligned}
\|\Delta_{-1}\tau(t)\|_{L^\infty}
&\lesssim\int_{|\xi|\le\frac{4}{3}}
|\widehat{\Delta_{-1}\tau}(t,\xi)|\,d\xi\\
&\lesssim(1+t)^{-\frac1\beta}
+\int_0^t(1+t-s)^{-\frac1\beta}
\left(\|F(s)\|_{L^1}+\|G(s)\|_{L^1}\right)\,ds\\
&\lesssim(1+t)^{-\frac1\beta}
\left(1+\int_0^{t/2}(1+s)^{-\frac{3}{2\beta}}\,ds\right)
+(1+t)^{-\frac{3}{2\beta}}
\int_{t/2}^t(1+t-s)^{-\frac1\beta}\,ds\\
&\lesssim(1+t)^{-\frac1\beta}.
\end{aligned}
\end{equation}
For $1/\beta<p<2$, Proposition \ref{3prop4} provides 
\begin{align}\label{3346}
\|G(s)\|_{L^p}+\|\nabla(u(s))\|_{L^2}\lesssim(1+s)^{-\frac1\beta}.
\end{align}
For the high-frequency part, applying
Lemma \ref{lem:heat-smoothing} and Bernstein's
inequality, together with \eqref{3237} and \eqref{3346}, we infer that
\begin{equation}\label{3347}
\begin{aligned}
\left\|\sum_{j\ge0}\dot\Delta_j\tau(t)\right\|_{L^\infty}
&\lesssim (1+t)^{-\frac1\beta}\|\tau_0\|_{\dot B^{-1}_{2,\infty}}
+\sum_{j\ge0}\int_0^te^{-c(t-s)2^{2\beta j}}
\left(2^{\frac{2j}{p}}\|G(s)\|_{L^p}
+2^j\|\nabla u(s)\|_{L^2}\right)\,ds\\
&\lesssim (1+t)^{-\frac1\beta}+\sum_{j\ge0}\left(2^{\frac{2j}{p}}+2^j\right)
\int_0^te^{-c(t-s)2^{2\beta j}}(1+s)^{-\frac1\beta}\,ds\\
&\lesssim (1+t)^{-\frac1\beta}
+\sum_{j\ge0}\left(2^{\frac{2j}{p}}+2^j\right)2^{-2\beta j}
\left(e^{-ct2^{2\beta j}/2}+(1+t)^{-\frac1\beta}\right)\\
&\lesssim(1+t)^{-\frac1\beta}.
\end{aligned}
\end{equation}
Combining \eqref{3345} and \eqref{3347}, we obtain
\[
\|\tau(t)\|_{L^\infty}
\lesssim (1+t)^{-\frac1\beta},
\qquad t\ge1.
\]
For $0\le t\le1$, the desired estimate follows
from the uniform $L^\infty$ bound \eqref{222}.
This completes the proof.
\end{proof}
\noindent\textbf{Proof of Theorem \ref{1theo2}:}
Theorem \ref{1theo2} follows from Propositions
\ref{3prop2}, \ref{3prop4} and \ref{3prop5}.
\hfill$\Box$

\smallskip
 
\noindent\textbf{Data Availability Statement.}
 Data sharing is not applicable to this article, as no data sets were generated or analyzed during the current study.

 \noindent\textbf {Conflict of Interest Statement} The authors declare that they have no conflict of interest. The authors also declare that this manuscript has not been previously published, and will not be submitted elsewhere before your decision.

\phantomsection
\addcontentsline{toc}{section}{\refname}
\bibliographystyle{abbrv} 
\bibliography{Oldroyd-Bref}

\end{document}